\documentclass[12pt]{amsart}

\usepackage{fullpage,url,amssymb,enumerate,colonequals}
\usepackage[all]{xy}
\usepackage{color}
\usepackage[OT2,T1]{fontenc}
\usepackage{graphicx}

\newcommand{\Ff}{\mathbb{F}_}

\newcommand{\Cc}{\mathbb{C}}
\newcommand{\Z}{\mathbb{Z}}

\newcommand{\Aa}{\mathbb{A}}
\newcommand{\mc}{\mathcal}
\newcommand{\mk}{\mathfrak}

\newcommand{\ov}{\overline}

\numberwithin{equation}{section}

\newtheorem{theorem}{Theorem}[section]
\newtheorem{lemma}[theorem]{Lemma}
\newtheorem{corollary}[theorem]{Corollary}

\theoremstyle{definition}

\newtheorem{examples}[theorem]{Examples}

\theoremstyle{remark}
\newtheorem{remark}[theorem]{Remark}

\definecolor{darkgreen}{rgb}{0,0.5,0}
\usepackage[
        colorlinks, citecolor=darkgreen,
        backref,
        pdfauthor={Zhiyao Zhang}, 
        pdftitle={Improving bounds for value sets of polynomials over finite fields},
]{hyperref}
\usepackage[alphabetic,backrefs,lite]{amsrefs} 

\begin{document}

\title{Improving bounds for value sets of polynomials over finite fields}

\author{Jiyou Li}

\address{School of Mathematical Sciences, Shanghai Jiao Tong University, 800 Dongchuan Road.}
\email{lijiyou@sjtu.edu.cn}

\author{Zhiyao Zhang}

\address{School of Mathematical Sciences, Shanghai Jiao Tong University, 800 Dongchuan Road.}
\email{aboctopus@sjtu.edu.cn}

\date{November, 2025}

\begin{abstract}
Let $\mathbb{F}_{q}$ be a finite field of characteristic $p$, and let $f \in \mathbb{F}_{q}[x]$ be a polynomial of degree $d > 0$. 
 Denote the image set of this polynomial as \(
V_{f}=\{f(\alpha)\mid\alpha\in\mathbb{F}_{q}\}
\)
and denote the cardinality of this set as $N_{f}$.
A much sharper bound for $N_{f}$ is established  in this paper. In particular, for any $p\neq 2, 3$, and for nearly every generic quartic polynomial  $f \in \Ff{q}[x]$, we obtain $$\lvert N_f - \frac{5}{8} q \rvert \leq \frac{1}{2}\sqrt{q} + \frac{33}{8},$$
which holds as a simple corollary of the main result.
\end{abstract}

\maketitle

\section{Introduction}

Let $\mathbb{F}_{q}$ be a finite field of order $q = p^{m}$ of characteristic $p$, and let $f \in \mathbb{F}_{q}[x]$ be a polynomial of degree $d > 0$. 
 Denote the image set of this polynomial as
\[
V_{f}=\{f(\alpha)\mid\alpha\in\mathbb{F}_{q}\}
\]
and denote the cardinality of this set as $N_{f}$.

The history of studying $N_{f}$ in various forms has a long history for over 160 years. 
Very few exact formulas for $N_{f}$ are known except for the polynomials of specific forms.  
Indeed,  computing $N_f$ for a given $f$ is computationally hard, with wide applications in number theory, coding theory and cryptography. 
For instance, even the special case of determining whether a polynomial is a permutation polynomial has drawn the attention of many mathematicians. 
 To date, no efficient algorithms are known for computing  $N_f$ exactly in general, and even no probabilistic polynomial-time algorithms are known \cite{cheng2013counting}.

Due to the hardness of exact computation, a major line of research seeks to bound $N_f$
   using algebraic information of  $f$.  For instance, the degree  $d$ of $f$ gives a trivial lower bound $N_f \geq \frac{q}{d}$. 
  An elegant and simple upper bound was given by Wan \cite{wan1993value}:
\[N_f\leq q-\frac{q-1}{d},\]
in the case that $f$ is not a permutation polynomial. In this paper we give a much tighter bound for quartic polynomials, which improves previous bounds significantly.

From now on we assume that $f(x)$ is separable, i.e. $f(x)$ is not a polynomial in $x^p$. Let $\Ff{q}(t)$ denote the rational function field in the indeterminate $t$, and let $G(f)$ (resp. $G^+(f)$) be the Galois group of the polynomial $f(x)-t$ over $\Ff{q}(t)$ (resp. $\ov{\Ff{q}}(t)$, where $\ov{\Ff{q}}$ is a fixed algebraic closure of $\Ff{q}$). A polynomial in $\Ff{q}[x]$ is called \textit{generic} if $G^+(f)$ is exactly the symmetric group $S_d$ on $d$ elements. In that case $G(f)=S_d$ holds  as well, because $G^+(f)$ is already the full permutation group.
Birch and Swinnerton-Dyer established the following significant result \cite{birch1959note}: for a generic polynomial $f$ of degree $d$ we have: $$N_f = \mu_dq+O(q^{1/2}),$$
where $\mu_d := \sum_{r=1}^d (-1)^{r-1}/r!$ and the implied constant in big O-notation is only dependent on $d$.

The result of Birch and Swinnerton-Dyer can be proved by applying the Riemann hypothesis for curves over finite fields to an irreducible component $C_r$ of the curve  $$Z_r: f(x_1) = f(x_2) = \ldots =f(x_r)$$ defined in the space $\mathbb{A}^{r}(\ov{\Ff{q}})$. One may instead use the Chebotarev density theorem for the corresponding
function-field extension, following an argument similar to that in
  \cite[pp. 145, Prop. 5.48]{iwaniec2021analytic}).

Experimental evidence suggests that the implied constant in \( O(q^{1/2}) \) should be relatively small compared to \( d \). 
Nevertheless, improving upon the trivial $O(d!)$ bound is generally challenging, and the optimal constant remains unknown, even for small $d$.  This is illustrated by the long-standing lack of improvement for the case  $d=4$. It's been long known that  
\[
N_f = \mu_3 q + O(1) = \frac{2}{3} q + O(1).
\]

But much less is known for \( d \geq 4 \). Using elementary methods, the authors of \cite{mccann1967distribution} proved the result of Birch and Swinnerton-Dyer for quartic polynomials without assuming the Riemann hypothesis for finite fields. In \cite{valentini2016note}, it was shown that
\[
N_f \geq \frac{q + 1}{2}
\]
for quartic polynomials of the form \( x^4 + a x^2 + b x \) with \( b \neq 0 \).

In the case $d=4$, The approach of Birch and Swinnerton-Dyer yields $$N_f = \frac{5}{8} q + O(q^{1/2}), $$
for a generic quartic $f$. After translating the variable and subtracting a constant, it suffices to consider depressed quartics of the form   $f(x) = x^4+ax^2 + bx$, with $b \neq 0$. We then define the normalized quantity: $$d_f =\frac{N_f - \frac{5}{8}q}{\sqrt{q}}.$$

\begin{figure}\label{Fig1}
    \centering
    \includegraphics[width=\linewidth]{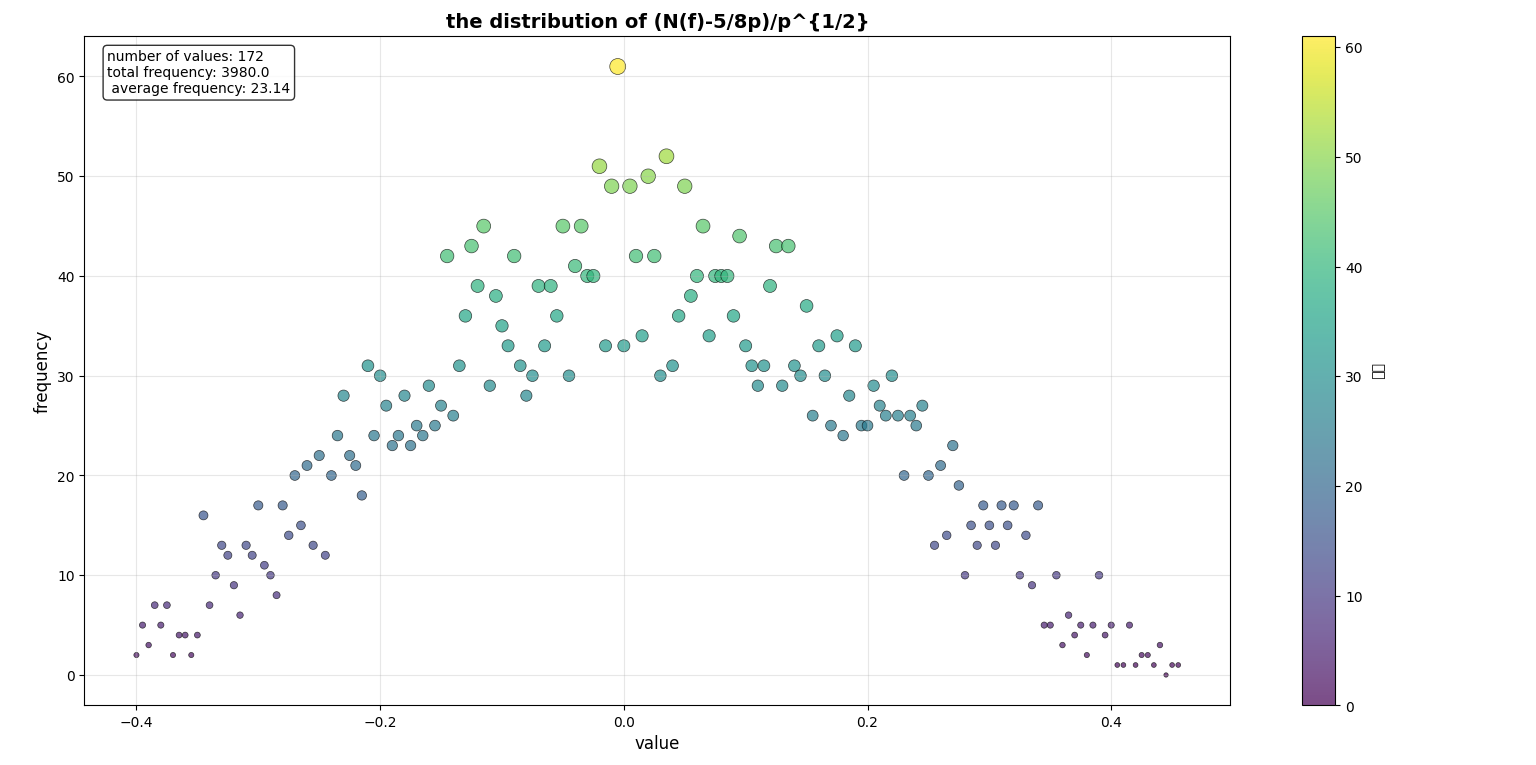}
    \caption{The distribution of implied constants}
    \label{fig:enter-label}
\end{figure}

To estimate the implied constant in  $O(q^{1/2})$, we examine experimental data. 
Using Magma, we compute the distribution of the values $d_f$ for the polynomial
  $f = x^4+ x^2 -x$ over various prime fields $\Ff{p}$, where  $p$ ranges over all primes $p \leq 37830$, as illustrated in FIGURE \ref{Fig1}.

The computational evidence suggests that the optimal square-root coefficient is likely around $1/2$. However, achieving this bound appears to be beyond the reach of existing methods.

 This paper is devoted to establishing a tighter bound for $d_f$  by refining the method of Birch and Swinnerton-Dyer. We establish a counting
identity, which expresses $N_f$ as an alternating linear combination of the numbers
$\#C_r(\Ff{q})$, where the $C_r$ are smooth projective models of some
fibre-product component of $Z_r$.  Artin formalism factors the zeta functions of the
$C_r$ into $L$-functions attached to representations of $S_d$.  Since the
same Artin factors occur for several values of $r$, the alternating
combination permits cancellation of their contributions.

The general theorem is intended as a framework for families in which the
component geometry and the relevant Artin conductors can be computed.  The
quartic case is treated separately because both calculations can be made
explicit there.  For degree five the conductor calculation becomes complicated and is not supplied here.

Our main result is the following theorem.

\begin{theorem} \label{mainthm2}
For a generic polynomial $f$ of degree $d$:
    $$\lvert N_f-\mu_dq \rvert \leq \sum_{\rho}\lvert \sum_{r=2}^d \frac{(-1)^r}{r!} m_{\rho,r} \rvert \deg(L(\rho,s)) \sqrt{q}  + C(d).$$
Here $\rho$ runs over all irreducible non-trivial representations of the symmetric group $S_d$, $m_{\rho,r}$ is the multiplicity of $\rho$ in the representation $\text{Ind}_{S_{d-r}}^{S_d} \chi_o$ induced from the trivial representation $\chi_o$, $L(\rho,s)$ is the Artin $L$-polynomial of $\rho:\text{Gal}(M/\Ff{q}(t)) \to GL(V)$ and $M$ is the splitting field of $f(x)-t$ over $\Ff{q}(t)$.
\end{theorem}

In the error term, the trivial representation is understood to contribute
zero: its zeta-factor supplies the main term $\mu_dq$, while the remaining
factors are the polynomial Artin $L$-functions described in Section 3.

Here $\Theta(d)$ denotes a quantity depending only on $d$ and on the fixed
degree bounds for the finite boundary and normalization terms. The explicit
restatement proved below takes $\Theta(d)=\frac{3}{2}d(e^d-d-1)+(\mu_d-1)$ before the additional
diagonal contribution is included.

The theorem is largely formal when it comes to polynomials of large degree. However, a concrete (albeit slightly trivial) bound can still be derived from our main theorem.

For clarity, a projective curve is called a complete intersection of
multidegree $(e_1,\ldots,e_{d-1})$ in $\mathbb P^d$ if it is the common zero
locus of a regular sequence of homogeneous equations of degrees
$e_1,\ldots,e_{d-1}$.  It is called smooth if it is nonsingular over
$\ov{\Ff{q}}$.

\begin{corollary} \label{concrete}
    The setup is the same as Theorem \ref{mainthm2}. If the projective closure of $T_d$ (See Section 3) in $\mathbb{P}^d$ is a smooth (ideal-theoretic) complete intersection, then: $$\lvert N_f - \mu_d q \rvert \leq 2^d (2+(d-1)!( \frac{d(d-1)}{2} -d-1 ))\sqrt{q}+\Theta(d).$$
\end{corollary}

The effectiveness of this theorem for small-degree polynomials is evidenced by the case $d=4$, where it leads to the following corollary:
 \begin{corollary}\label{quar} Suppose $p\neq 2, 3$, there exists a Zariski open dense subset (see Section 4) $U \subseteq \Aa^2$ such that: for any $(a,b) \in U$, if $f=x^4+ax^2+bx \in \Ff{q}[x]$ is generic, we have $$\lvert N_f - \frac{5}{8} q \rvert \leq \frac{1}{2}\sqrt{q} + \frac{33}{8}.$$
\end{corollary}

Note that our result gives that asymptotically $\lvert d_f \rvert \leq \frac{1}{2}$.  This bound holds for nearly all generic quartic polynomials, fits the experimental data very well and improves upon all previous results in the literature.

The condition char $\Ff{q} \neq 2,3$ is necessary. When char $\Ff{q}=2$, $f$ degenerates into a linear polynomial, and its image set is a $\Ff{p}$-vector space.

This paper is structured as follows. Section 2 reviews the method of Birch and Swinnerton-Dyer. Section 3 covers the prerequisite concepts and foundational theorems of Artin L-functions. Section 4 is devoted to the proof of Theorem 1.1 (= Theorem \ref{mainthm}) and proves Corollary \ref{concrete}. Finally, Section 5 presents the proof of Corollary  \ref{quar}. 

\section{The method of Birch and Swinnerton-Dyer}

Let notations and assumptions be as in the previous section, and the polynomial $f$ will always be a monic and separable one of degree $d$. We place no restriction on the Galois group of $f(x)-t$ over $\Ff{q}(t)$ until the generic case is considered. In this section, we briefly recap the theorem by Birch and Swinnerton-Dyer and its proof. Compared to the original paper of Birch and Swinnerton-Dyer, the method employed in our paper is more profound and reveals the decomposition structure of the curves with greater clarity. Furthermore, our approach extends a commonly used method of studying polynomial properties, namely the \textit{monodromy} method \cite{fried1993schur}.

\subsection{An initial estimate of $N_f$} This subsection is a paraphrase of Birch and Swinnerton-Dyer's original paper. For every integer $r$ with $2 \leq r \leq d$, consider the curves $Z_{r}$ defined in $\mathbb{A}^{r}(k)$ ($k = \ov{\Ff{q}}$): $$f(x_1) = f(x_2) = \ldots = f(x_{r}).$$ Regard $f$ as a map from affine $x$ line to affine $t$ line $f: \mathbb{A}_x^1 \to \mathbb{A}_t^1, x \mapsto t=f(x)$. By definition, $Z_r$ is the fibre product of $\mathbb{A}_x^1$ by itself $r$ times: $$Z_r = \underbrace{\mathbb{A}_x^1 \times_{\mathbb{A}^1_t} \times \mathbb{A}_x^1 \times_{\mathbb{A}^1_t} \ldots  \times_{\mathbb{A}^1_t} \mathbb{A}_x^1}_{r \text{ times}}. $$ It is clear that the coordinate ring of $Z_r$ has the isomorphism \begin{align*}k[x_1,x_2,\ldots,x_r,t]/(t-f(x_1),t-f(x_2),\ldots,t-f(x_r)) \\ \cong k[x_1] \otimes_{k[t]} k[x_2] \otimes_{k[t]} \ldots k[x_r] .\end{align*}

$Z_r$ decomposes into a finite union of irreducible curves defined over $k$. For example, the curve $Z_2$ decomposes into at least two components, $x_1 = x_2$ and $$\frac{f(x_1)-f(x_2)}{x_1-x_2} = 0.$$ Any component is defined over $\Ff{q}$ iff. it is closed under the Frobenius polynomial mapping $x_i \mapsto x_i^q$. Let $T_r$ be the union of irreducible components of $Z_r$ not completely contained in any of the hyperplane $x_i = x_j$. $T_r$ is closed under Frobenius mapping, therefore $T_r$ is defined over $\Ff{q}$. However, its irreducible components are not necessarily defined over $\Ff{q}$. Let $n_r$ denote the number of $\Ff{q}$-rational points on $T_r$, and let $n_r'$ denote the number of $\Ff{q}$-solutions to $$f(x_1) = f(x_2) = \ldots = f(x_r)$$ such that $x_1, \ldots,x_r$ are distinct. The projective closure of $T_r$ has degree at most $d^{r-1}$. Since no component of $T_r$ is contained in any hyperplane $D_{ij}=\{x_i=x_j\}$, B\'ezout's theorem gives $$\#(T_r \cap D_{ij})(\Ff{q}) \leq d^{r-1}.$$ There are $\binom{r}{2}$ hyperplanes defined by $x_i=x_j$, so consequently $n_r = n_r'+O_d(1)$, with an explicit bound depending only on $d$.

Denote by $m_i$ the number of $y$ that has $i$ distinct roots of $f(x)=y$ over $\Ff{q}$. Certainly  $$N_f = m_1 + m_2 + \ldots +m_d;$$ $$q = m_1 + 2m_2 + \ldots + dm_d.$$
One can easily count $n_r'$ in terms of $m_i$. $$n_r' = r! m_r + \frac{(r+1)!}{1!}m_{r+1} + \ldots + \frac{d!}{(d-r)!}m_{d}.$$
And thus 
$$\frac{n_r'}{r!} = m_r + \binom{r+1}{r} m_{r+1} + \ldots + \binom{d}{r}m_{d}.$$

So we have  \begin{align*} \sum_{r=2}^d (-1)^{r} \frac{n_r'}{r!} &= (-1+2) m_2 + \ldots + (\sum_{k=0}^{d-2} (-1)^{d-k}\binom{d}{k}) m_d \\ & = m_2 + 2m_3 + \ldots + (d-1)m_d \\ & = q - N_f.
\end{align*}

On the other hand, by Hasse-Weil bound, we have: $$n_r = \nu_r q + O(q^{1/2}).$$ where $\nu_r$ is the number of irreducible components of $T_r$ over $k$ defined over $\Ff{q}$, and the implied constant in $O(q^{1/2})$ is dependent on the geometric genus of irreducible components of $T_r$. It is well known that those components not defined over $\Ff{q}$ will at most have contribution $O_d(1)$ to $n_r$. (Suppose $C$ is defined over $k$, but not over $\Ff{q}$, and $C$ is cut out by polynomials $f_1,f_2,\ldots,f_r$. Let $\ov{C}$ be the \textit{Frobenius twist} of $C$. It is cut out by $\ov{f_1}, \ov{f_2},\ldots,\ov{f_r}$, where each $\ov{f_i}$ is obtained from $f_i$ by replacing each coefficient of $f_i$ with its $q$-th power. Then any $\Ff{q}$-point of $C$ lies in the intersection $C \cap \ov{C}$, whose cardinality $\lvert C \cap \ov{C} \rvert = O_d(1)$.)

The difference between $q-N_f$ and $\sum_{r=2}^d (-1)^r \frac{n_r}{r!}$ is bounded by:
\begin{align*} \lvert \sum_{r=2}^d (-1)^{r} \frac{n_r'-n_r}{r!} \rvert & \leq \lvert\sum_{r=2}^d \frac{d^{r-1} r(r-1)}{2r!} \rvert \\ & = \frac{1}{2} \sum_{r=2}^{d} \frac{d^{r-1}}{(r-2)!} \\ & \leq \frac{d}{2} (e^d-d-1).  \end{align*}
This shows that the error is dependent on $d$.
Combined with Hasse-Weil bound, we establish a general formula for the cardinality of the image set of an arbitrary $f$ $$N_f = (1+\sum_{r=2}^d (-1)^{r+1}\frac{\nu_r}{r!}) q + O(q^{1/2}).$$
The value $\nu_r$ is controlled by the Galois group $G^+(f)$ of $f(x)-t$ over $k(t)$ and $G(f)$ (Galois group over $\Ff{q}(t)$). Generally, let $M$ be the splitting field of $f(x)-t$ over $\Ff{q}(t)$, and denote by $\ell$ the algebraic closure of $\Ff{q}$ inside $M$. By definition, $G^+(f)$ is the normal subgroup of $G(f)$ leaving the subfield $\ell \subseteq M$ fixed, and we have an isomorphism $G(f)/G^+(f) \cong Gal(\ell/\Ff{q})$.

Both $G(f)$ and $G^+(f)$ act naturally on the set of roots $X:=\{x_1,x_2,\ldots,x_d\}$ of $f(x)-t$ transitively. Denote by $X^m$ the set of $m$-tuples of distinct elements of $X$, on which $G(f)$ and $G^+(f)$ also act naturally.

In fact, $\nu_r$ counts the number of $G(f)$-orbits of $X^r$ that are also $G^+(f)$-orbits. Their argument is as follows: On the curve $\{f(x_1)=\ldots = f(x_r)=t\} \cong Z_r$, let $P=(\zeta_1,\zeta_2,\ldots, \zeta_r,\eta)$ be such that $\eta$ is \textit{generic} over $k$ and $\zeta_1, \ldots,\zeta_r$ are distinct roots of $f(x)=\eta$. Then $P$ is a simple point on $T_r$, and a unique component $C$ passes through $P$. The group $G^+(f)$ and $G(f)$ acts naturally on $P$, and $C$ contains precisely $\sigma^+P$ for $\sigma^+ \in G^+(f)$. If each $\sigma P$ such that $\sigma \in G(f)$ is also a $\sigma^+P$ for some $\sigma^+ \in G^+(f)$, then $C$ is defined over $\Ff{q}$. Interested readers can refer to Birch \& Swinnerton-Dyer's paper for the details.

\begin{examples}
    When $f$ is a \textit{generic} polynomial (i.e. $G^+(f) = S_d$), clearly $G(f)$ is also $S_d$. For every value of $r$, $G^+(f)$ has only one orbit on $X^r$. Therefore, $T_r$ is irreducible and defined over $\Ff{q}$, and $\nu_r = 1$ for $2 \leq r \leq d$.
\end{examples}

\begin{examples}
    Let $f(x)=x^d$ such that $(d,q-1)=1$. $$f(x)-f(y)=x^d-y^d = \prod_{r=0}^{d-1} (x-\mu_d^r y).$$ breaks into linear polynomials over $k[x,y]$, where $\mu_d$ is a primitive $d$-th root of unity of $\Ff{q}$. Since $(d,q-1)=1$, $x - \mu_d^r y$ is not defined over $\Ff{q}$ except when $r=0$. By definition, $T_2$ contains the lines $x - \mu_d^r y = 0$ defined over $k$ ($r \geq 1$), but none of them are defined over $\Ff{q}$, so $\nu_2=0$. In this particular case, $G^+(f) = \Z/d\Z$ and $G(f) = \Z/d\Z \rtimes (\Z/d\Z)^*$.
    
    $f$ is an example of exceptional polynomials over $\Ff{q}$: every $G(f)$-orbit of $X^2$ breaks into smaller $G^+(f)$-orbits.
\end{examples}

\begin{remark}
The method of Birch and Swinnerton-Dyer can be generalized to a broader setting. Suppose we have a function $\omega:\Ff{q} \to \Cc$ where $$V=\sum_{x \in \Ff{q}}\omega(f(x))$$ is relatively easier to compute or estimate. For example, $\omega$ can be $\chi:(\Ff{q},+) \to \Cc^*$ a non-trivial additive character. In this case, there is a non-trivial estimate $|V| \leq (d-1)\sqrt{q}$ by Weil's bound. Suppose we want to estimate $$V' = \sum_{y \in V(f)} \omega(y).$$

For $2 \leq r \leq d$, consider the following values  $$V_r =\sum_{x \in N_r} \omega(f(x)).$$ $N_r$ is the set of $\Ff{q}$-rational points on $T_r$ with distinct coordinates, and here $f(x)$ should be understood as a function $N_r \to \Ff{q}: (x_1,x_2,\ldots,x_r) \mapsto f(x_1)$. The completely same argument as before shows $$V - V' = \sum_{r=2}^d (-1)^{r} V_r/r!.$$

In the case where $\omega$ is a non-trivial additive character on $\Ff{q}$, we have Weil's bound $$V \leq (d-1)\sqrt{q}.$$ and $V_r$ is an exponential sum on a curve.
\end{remark}

\subsection{The decomposition structure of $Z_r$}
In this subsection, we explain how the curve $Z_r$ breaks into a lot of irreducible components from a ring theoretic viewpoint, so that the readers can understand the structure of $Z_r$ and $T_r$ more clearly. One major takeaway of this subsection is that when $f$ is generic, the function field of $T_r$ is $\Ff{q}(x_1,x_2,\ldots,x_r)$, where $x_1,x_2,\ldots,x_r$ are distinct roots of $f(x)-t$.

\begin{lemma}
The irreducible components of $Z_r$ over $k$ corresponds to the direct components of the tensor product of fields: $$k(x) \otimes_{k(t)} k(x) \ldots \otimes_{k(t)} k(x).$$ The correspondence is $$V \mapsto \text{Function field of }V.$$
\end{lemma}
\begin{proof}
Since $\Aa_x^1 \to \Aa_t^1$ is finite flat, $Z_r \to \Aa_t^1$ is also flat, and the generic point $\xi$ of any irreducible component $V$ of $Z_r$ maps to the generic point $\eta$ of $\Aa^1_t$. \cite[Part 2, Lemma 29.26.9]{stacks-project} Conversely, for any point $\xi$ in the generic fibre of $Z_r \to \Aa_t^1$, its Zariski closure in $Z_r$ is an irreducible curve in $Z_r$, hence an irreducible component of $Z_r$. The generic fibre of $Z_r \to \Aa_t^1$ is isomorphic to $Z_r \times_{\Aa_t^1} \text{Spec }k(t) \cong \text{Spec}(k(x) \otimes_{k(t)} k(x) \ldots \otimes_{k(t)} k(x))$. This proves the first part of lemma.

It is well known that tensor products of fields are direct sum of various fields (See for example, \cite[Section 8.18]{Jacobson1989}). By definition of generic points, these fields are function fields of the irreducible components of $Z_r$. 
\end{proof}

It's worth pointing out that following the theory of tensor products of fields, the function field of any irreducible component of $Z_r$ looks like
$$k(x_1,x_2,\ldots,x_m)$$ where $x_i$ are distinct roots of $f(x)-t$ and $m \leq r$.

Let us look at the irreducible components of $Z_r$ in more detail. First, we look at $Z_2:f(x_1)=f(x_2)$ as an example.
Suppose $Z_2$ is the (scheme-theoretic) union of irreducible components $x_1=x_2$ and $f_i(x_1,x_2)=0$. Every irreducible component of $Z_2$ (except $x_1=x_2$ itself) is not completely contained in the subspace $x_1=x_2$. Furthermore, the coordinate rings of $f_i(x_1,x_2)=0$ are generated by two distinct roots $x_1$ and $x_2$ of $f(x)-t$, considering $f(x_1)=f(x_2)=t$.

Let's figure out the decomposition of $Z_3:f(x_1)=f(x_2)=f(x_3)=t$. We should consider the components $$\{f_i(x_1,x_2)=0\} \times_{\Aa_t^1} \Aa_x^1.$$ and $\{x_1=x_2\} \times_{\Aa_t^1} \Aa_x^1$ separately. The latter is simple: it is isomorphic to $Z_2$ and completely contained in $D_{12}:x_1=x_2$. It stands for the component $$f(x_1)=f(x_2)=f(x_3),x_1=x_2.$$ The coordinate ring of the former component is: 
$$
    k[x_1,x_2]/(f_i(x_1,x_2))) \otimes_{k[t]} k[x_3] \\ 
     \cong R[x_3]/(f(x_3)-t).
$$ by writing $R=k[x_1,x_2]/f_i(x_1,x_2)$ a domain and $t=f(x_1)=f(x_2)$.
We decompose the polynomial $f(x)-t$ with indeterminate $x$ over $R$: Suppose the decomposition is $(x-x_1)(x-x_2) \prod_j g_j(x_1,x_2,x)$, then the irreducible components of Spec $k[x_1,x_2] \otimes_{k[t]}k[x_3]$ are given by $$f_i(x_1,x_2)=0, x_j =x_3 (j=1,2), \text{ coordinate ring }k[x_1,x_2]$$ and $$f_i(x_1,x_2)=0, g_j(x_1,x_2,x_3)=0, \text{ coordinate ring }k[x_1,x_2,x_3].$$ 

More generally, $Z_r$ is the union of various curves whose coordinate rings are obtained from $k[t]$ by adjoining specific roots of $f(x)-t$. As an affine variety, each component can be defined in $\Aa^r$ using $r-1$ equations, each having degree $d-1,d-2,\ldots,d-r+1$. (The reader is recommended to read Section 4 for a concrete example.) Each direct summand stands for an irreducible component of $Z_r$ over $k$. The number of irreducible components of $Z_r$ is strongly reliant on the behavior of irreducible decomposition of $f(x)-t$ over various subrings of $k[x_1,\ldots,x_d]$, where $x_1,\ldots,x_d$ are certain roots of $f(x)-t$ over $k(t)$.

An irreducible component of $Z_r$ is contained in some $x_i = x_j$ if and only if its coordinate ring can be obtained by adding less than $r$ roots of $f(x)-t$ to $k$.

\begin{examples}
    When $f$ is generic, we have $G(f)=G^+(f)=S_d$. Therefore, the groups $G(f)$ and $G^+(f)$ acts transitively on the sets $X^r$. Since the $G^+(f)$-orbits (resp. $G(f)$-orbits) count the number of irreducible components of $T_r$ over $k$ (resp. over $\Ff{q}$), it immediately follows that $T_r$ is irreducible and defined over $\Ff{q}$. Following a similar analysis to the case over $k$, the coordinate ring of $T_r$ (over $\Ff{q}$) is $\Ff{q}[x_1,x_2,\ldots,x_r]$, where $x_i$ are distinct roots of $f(x)-t$ over $\Ff{q}(t)$.
\end{examples}

\section{Preliminaries from algebraic number theory}

In this section, we recall necessary concepts and theorems from algebraic number theory in global function fields, i.e. finite extension of $\Ff{q}(t)$. \cite{rosen2013number}

First we recall the concept of a Frobenius conjugacy class. For a Galois extension $L/K$ of global function fields with Galois group $G$, $P$ a prime of $K$, and $\mk{P}$ a prime of $L$ lying above $P$, their residue fields $K_P \hookrightarrow L_\mk{P}$ are finite, and the Galois group Gal$(L_\mk{P}/K_P)$ is cyclic of degree $f(\mk{P}/P)$, which is the inertia degree of $\mk{P}$ over $P$. It is generated by the Frobenius automorphism $\phi_P:x \mapsto x^{|K_P|}$.

Define the decomposition group $D(\mk{P}/P)$ as the subset of $G$ that fixes $\mk{P}$, and the inertia group $I(\mk{P}/P)$ as the subset of $D(\mk{P}/P)$ that induces identity on the residue class field $L_\mk{P}$, i.e. $\sigma \in G$ such that $\sigma \cdot x \equiv x \mod \mk{P}$ in the integer ring of $L$. We have the exact sequence of groups $$1 \to I(\mk{P}/P) \to D(\mk{P}/P) \to \text{Gal}(L_\mk{P}/K_P) \to 1.$$
Specifically if $P$ is unramified in $L$, then $I(\mk{P}/P)$ is trivial, and $D(\mk{P}/P) \to \text{Gal}(L_\mk{P}/K_P)$ becomes an isomorphism. The pre-image of $\phi_P$ under this isomorphism is called the Frobenius automorphism of $\mk{P}$ for the extension $L/K$, denoted $(\mk{P},L/K)$. For various primes $\mk{P}$ above $P$, $(\mk{P},L/K)$ fills a conjugacy class in $G$, called the Frobenius conjugacy class of $P$, denoted Frob$_P$.

Fix a finite dimensional linear representation $\rho:G \to GL(V)$ over $\Cc$. For $P$ a prime of $K$ unramified in $L$, define the local factor at $P$.
$$L_P(\rho,s):= \det(I-\rho(\text{Frob}_P) N(P)^{-s})^{-1}.$$ Here $N(P)$ denotes the norm of prime ideals. For $P$ ramified in $M$, fix a prime $\mathfrak{P}$ above $P$, and denote the invariant subspace of $V$ under the inertia subgroup $I(\mathfrak{P}/P)$ by $V^I$. Instead of defining the local factors as above we define: $$L_P(\rho,s):= \det(I-\rho(\gamma_P)\rvert_{V^I} N(P)^{-s})^{-1}.$$ $\gamma_P$ is any element in $D(\mathfrak{P}/P)$ mapping to the Frobenius element in Gal$(L_{\mathfrak{P}}/K_P)$. This definition is independent of the element chosen.

Finally, the Artin L-function $L(\rho,s)$ is defined as $$L(\rho,s):= \prod_P L_P(\rho,s).$$
where $P$ ranges over all primes of $K$.

The basic properties of Artin L-functions we will use are as follows. They can all be found in \cite[Chapter 9]{rosen2013number} except the last two. The last two can be found on \cite[Page 4]{}, for example.

\begin{itemize}
    \item For the trivial representation $\rho = \chi_o$, we have $L(\chi_o,s) =\zeta_K(s)$, the zeta function of $K/\Ff{q}$. Riemann's hypothesis over finite fields shows that there exists a polynomial $P_K$ in $u=q^{-s}$ such that: $$\zeta_K(s) = \frac{P_K(u)}{(1-u)(1-qu)}.$$ Here $$P_K(u) = 1+ a_1u+a_2u+\ldots + a_{2g}u^{2g} = \prod_{i=1}^{2g}(1-\alpha_i u).$$ splits into linear factors over $\Cc$. The degree $g$ of $P_K$ is the genus of $K$, and the complex modulus of every root of $P_K(u)$ is $q^{-1/2}$. In other words, $\lvert \alpha_i \rvert_{\Cc}=q^{1/2}$ for every $i$.
    \item For the regular representation $\rho = \mathbb{C}[G]$, we have $L(\rho,s) = \zeta_L(s)$.
    \item For two representations $\rho_1$ and $\rho_2$, we have $L(\rho_1 \oplus \rho_2,s) = L(\rho_1,s) L(\rho_2,s)$. To study the $L$-function of any representation, it suffices to study the $L$-functions of irreducible representations.
    \item Given a subgroup $H$ of $G$, denote the fixed field of $H$ in $L$ by $L^H$. Let $\sigma$ be a $\Cc$-representation of $H = \text{Gal}(L/L^H)$, and $\text{Ind}_H^G \sigma$ be the induced representation. Then $L_G(\text{Ind}_H^G \sigma,s) = L_{H}(\sigma,s)$.
    \item In the case of $\sigma = \chi_o$ the trivial representation, we have $L_G(\text{Ind}_H^G \chi_o,s) = \zeta_{L^H}(s)$.
\end{itemize}

The regular representation $\Cc[G]$ decomposes into several irreducible representations $\bigoplus_\rho \dim(\rho) \rho$, where $\rho$ ranges over all irreducible representations of $G$. We have the equality $$\zeta_L(s) = \zeta_K(s) \prod L(\rho,s)^{\dim \rho}.$$ Multiplying both sides by $(1-u)(1-qu)$ $$P_L(u) = P_K(u) \prod L(\rho,s)^{\dim \rho}.$$

It's well known that if $L/K$ is a geometric extension, then $L(\rho,s)$ is a polynomial in $u =q^{-s}$ for non-trivial irreducible representations. Denote by $Z(u)$ the polynomial satisfying $Z(q^{-s})=L(\rho,s)$. By comparing two sides of the equation, we easily know: 
\begin{itemize}
    \item The roots of polynomial $Z(u)$ over $\Cc$ also have modulus $q^{-1/2}$ (if there is any).
\end{itemize}

The degree of $Z(u)=L(\rho,s)$ is given by Raynaud’s Euler characteristic formula \cite{raynaud1964caracteristique}
$$\deg L(\rho,s) = (2g-2) \dim \rho + \deg f(\rho).$$ where $f(\rho)$ is the degree of the \textit{Artin conductor} of $\rho$. When $g = 0$ (rational function field $K$), $f(\rho) \geq 2 \dim \rho$, and the equality holds iff. $L(\rho,s) = 1$ the constant polynomial.

To conclude this section, we recall the concept of Hasse-Weil zeta functions of an algebraic curve (\cite[Chapter 5]{rosen2013number} and \cite{niederreiter2001rational}). Let $X$ be an algebraic curve defined over $\Ff{q}$, denote by $\#X(\Ff{q^r})$ the cardinality of the set of its $\Ff{q^r}$-rational points. Define its Hasse-Weil zeta function $$\zeta(X,s):= \exp(\sum_{n=1}^\infty \frac{\#X(\Ff{q^n})}{n}q^{-sn}).$$ By Weil's theorem, $\zeta(X,s)$ is a rational function in $q^{-s}$. Writing $T=q^{-s}$, we have $$\zeta(X,s) = \frac{P_X(T)}{(1-T)(1-qT)}.$$ where $P_X(T)$ is a polynomial with integer coefficients.
 Denote by $K$ the function field of $X$. When $X$ is projective and smooth, we have the equality $\zeta(X,s) = \zeta_K(s)$. Conversely, every function field $K$ of one variable over $\Ff{q}$ is uniquely associated up to isomorphism to a smooth projective curve $X$ whose function field is isomorphic to $K$; call it the smooth projective model of $K$.

\section{Proof of Theorem \ref{mainthm}}

In this section, we combine the material developed in the preceding sections and develop an expanded analytical framework of image sets of a generic polynomial $f$. The notations are the same as Section 2, and we always suppose $f$ is a generic polynomial. When $f$ is generic, $T_r$ is absolutely irreducible, defined over $\Ff{q}$, and has function field $K_r = \Ff{q}(x_1,\ldots,x_r)$, where $x_1,x_2,\ldots,x_r$ are distinct roots of $f(x)-t$ over $\Ff{q}(t)$. There is a tower of function fields $$\Ff{q}(t) \subseteq K_1 \subseteq K_2 \subseteq \ldots \subseteq K_d =M.$$ 
Here $M/\Ff{q}(t)$ is a Galois extension with Galois group $G =S_d$. The degree of each extension is $$[K_i:\Ff{q}(t)] = d (d-1) \ldots (d-i+1), i \leq d-1.$$ Moreover, $K_{d-1}=K_d$. The smooth projective model of $K_r/\Ff{q}$ is birational to $T_r$, since they share the same function field.

As a subfield of $M$, $K_r (1 \leq r \leq d-1)$ is the fixed field of $S_{d-r} \subseteq G = S_d$, where $S_{d-r}$ can be identified with the automorphisms only permuting the roots $x_{r+1}, \ldots,x_d$. Suppose the induced representation $\text{Ind}_{S_{d-r}}^{S_d}\chi_o$ decomposes into irreducible representations $\bigoplus_\rho m_{\rho,r}\rho$, where $m_{\rho,r}$ is the multiplicity of $\rho$. By the induced representation property of Artin $L$-functions (introduced in the last section) $$\zeta_{K_r}(s) = L(\text{Ind}_{S_{d-r}}^{S_d} \chi_o,s) = \frac{\prod_\rho L(\rho,s)^{m_{\rho,r}}}{(1-u)(1-qu)}.$$
 $\rho$ runs over all non-trivial representations. We omit trivial representations since $L(\chi_o,s) = 1$ is trivial. Since $K_{d-1} = K_d$, we know that $m_{\rho,d} = m_{\rho,d-1} = \dim \rho$.

For our convenience, we introduce a notation: for a polynomial $P(t) = \prod_{i=1}^d (1-\omega_i t)$ and a positive integer $n$, define $$U(P,n):=\sum_{i=1}^d \omega_i^n.$$ Then the relation between the number of $\Ff{q^n}$-rational points on a smooth projective curve$/\Ff{q}$ and its Hasse-Weil zeta function can be described as \begin{equation} \label{weil}
    \zeta(X,s) = \frac{P_X(u)}{(1-u)(1-qu)} \Longrightarrow \#X(\Ff{q^n}) = q^n + 1 -U(P_X,n). \tag{*}
\end{equation} Call $P_X$ the $L$-polynomial of $X$.
From the formula $$N_f = q + \sum_{r=2}^d (-1)^{r+1} \frac{n_r'}{r!}.$$
The smooth projective model $C_r$ of $K_r/\Ff{q}$ differs from $T_r$ by a finite set of closed points, coming from singularities of $T_r$ and points of $C_r$ at infinity. Since the maximal number of singularities in $T_r$ is bounded by its arithmetic genus, which is less than $(d-1)d^r$, and the curve has degree $\leq d^r$, so the difference between $C_r$ and $T_r$ comes from $d^{r+1}$ points at most. So the total difference between $\sum_{r=2}^d(-1)^{r+1} \frac{n_r}{r!}$ and $\sum_{r=2}^d(-1)^{r+1} \frac{\#C_r(\Ff{q})}{r!}$ is bounded by  $$\sum_{r=2}^d d^{r+1}/r! \leq C(d) = d(e^{d}-d-1).$$ 
Since $C(d)$ is very large, the total difference could be a large value. But for most curves, this value is expected to be very small since $T_r$ is likely to be smooth or contain fewer singularities. As we have explained in Section 2, the difference between $\sum_{r=2}^d(-1)^{r+1} \frac{n_r}{r!}$ and $q-N_f$ is bounded by $1/2 \cdot C(d)$. Hence by substituting $n_r' = \#C_r(\Ff{q})+O_d(1)$, we are led to Theorem \ref{mainthm}:
\begin{align*}
N_f &= q + \sum_{r=2}^d (-1)^{r+1} \frac{\#C_r(\Ff{q})}{r!} + \frac{3}{2}d(e^d-d-1) \\ &(\text{Let }X = C_r \text{ in \ref{weil}}) \\ &= \mu_dq + \sum_{r=2}^d \frac{(-1)^r}{r!} U(P_{C_r}, 1) + \frac{3}{2} d(e^d-d-1)+\mu_d-1 \\ & = \mu_dq + \sum_{r=2}^d \frac{(-1)^r}{r!} \sum_\rho m_{\rho,r} \cdot U(L(\rho,s), 1) + \frac{3}{2} d(e^d-d-1)+\mu_d-1 \\ &= \mu_dq + \sum_{\rho} (\sum_{r=2}^d \frac{(-1)^r}{r!}  m_{\rho,r}) \cdot U(L(\rho,s), 1) + \frac{3}{2} d(e^d-d-1)+\mu_d-1.
\end{align*}

From our preceding discussion of Artin $L$-functions, $U(L(\rho,s),1)$ is the sum of $\deg L(\rho,s)$ algebraic integers of complex modulus $\sqrt{q}$.

\begin{theorem}\label{mainthm} For $f$ a generic polynomial of degree $d$:
 $$\lvert N_f-\mu_dq \rvert \leq \sum_{\rho}\lvert \sum_{r=2}^d \frac{(-1)^r}{r!} m_{\rho,r} \rvert \deg(L(\rho,s)) \sqrt{q} + \frac{3}{2}d(e^d-d-1)+(\mu_d-1).$$
\end{theorem}

\begin{proof}[Proof of Corollary \ref{concrete}]
Firstly, we have:
    \begin{align*}
        \lvert \sum_{r=2}^d \frac{(-1)^r}{r!} m_{\rho,r} \rvert \leq \sum_{r=2}^d \frac{m_{\rho,r}}{r!} & \leq \sum_{r=2}^d \frac{\dim P_r/\dim \rho}{r!} \\ & \leq \frac{1}{\dim \rho} \sum_{r=2}^d \frac{d!}{(d-r)!r!} \leq \frac{1}{\dim \rho} 2^d \leq 2^d.
    \end{align*} Therefore: \begin{align*}
        \lvert N_f-\mu_dq \rvert \leq 2^d \cdot \sum_{\rho} \deg L(\rho,s) \sqrt{q}+ C(d) \leq 2^d (2+(d-1)!( \frac{d(d-1)}{2} -d-1 ))\sqrt{q}+C(d).
    \end{align*} This is because $\sum_\rho \deg L(\rho,s) \leq \sum_\rho \dim \rho \cdot\deg L(\rho,s) \leq \deg P_{K_d}$. The projective closure of $T_d$ in $\mathbb{P}^d$ is smooth, so it is a complete intersection of multidegree $(d-1,d-2,\ldots,1)$ isomorphic to $C_d$, whose genus is well known to be: $$g =1+\frac{1}{2} (d-1)! (\sum_{i=1}^{d-1} i - d-1) = 1+ \frac{1}{2}(d-1)!( \frac{d(d-1)}{2} -d-1 ). $$ Hence $\deg P_{K_d} = 2g=2+(d-1)!( \frac{d(d-1)}{2} -d-1 )$.
\end{proof}

The projective closure of $T_d$ is smooth for almost all generic $f$. It is a "set theoretic complete intersection" in $\mathbb{P}^d$, i.e. can be cut out by $d-1$ equations $F_1,F_2\ldots,F_{d-1}$ (with variable $x_0,x_1,\ldots,x_d$). When the projective closure of $T_d$ is smooth, the polynomial ideal generated by $F_1,F_2,\ldots,F_{d-1}$ in $\Ff{q}[x_0,x_1,\ldots,x_d]$ is well-known to be radical, and the condition of Corollary \ref{concrete} holds.

This estimate concerns the sum over the irreducible representations of $S_d$, which correspond one-to-one to partitions of $d$. In the right hand side of the equation, there seems to be too many terms to carry out any meaningful estimate. Thankfully, some of those terms may be zero. There are two cases, either $m_{\rho,r}$ is zero or $\deg(L(\rho,s))$ is zero.

\subsection{Combinatorial interpretation of $m_{\rho,r}$}

As a $S_n$-module, the induced representation $P_r :=\text{Ind}_{S_{d-r}}^{S_d} \chi_o$ of the trivial representation is isomorphic to the free $\Cc$-vector space generated by ordered $r$-tuples of distinct elements in $\{1,2,\ldots,n\}$.

The modules $P_r$ are so-called $(d-r,1^r)$-\textit{Young permutation modules} in the representation theory of $S_n$ \cite[Chapter I.4]{fulton2013representation}. For a partition $\lambda$ of $n$, denote by $\Delta^\lambda$ the set of all Young tableau of shape $\lambda$. Define an equivalence relation $\sim$ on $\Delta^\lambda$: two tableaux are equivalent iff. their rows are the same. Finally, define $M^\lambda$ to be the free $\Cc$-vector space generated by $\Delta^\lambda/\sim$ with natural $S_n$-action.

It is well known that the irreducible representations of $S_n$, i.e. Specht modules can also be indexed by partitions $\mu$ of $n$, denoted $S_\mu$ for example. $M^\lambda$ decomposes into some irreducible modules $$M^\lambda = \bigoplus K_{\mu \lambda} S_\mu.$$

$K_{\mu \lambda}$ is the Kostka number of $(\mu, \lambda)$, equal to the number of semistandard Young tableaux (SSYT) of shape $\mu$ and weight $\lambda$.

To sum up, $m_{S_\mu,r} = K_{\mu\lambda}$ is a Kostka number where $\lambda =(r,1,\ldots,1)$. To put more detail to it, it counts the number of fillings of $r$ one's, and $2, \ldots, n-r+1$ into the Young diagram of shape $\mu$ such that the entries weakly increase along each row and strictly increase down each column. 

If $\lambda = (\lambda_1, \ldots,\lambda_r)$, and $\lambda_1 < r$, then it is clear $K_{\lambda \mu} = 0$ since the $r$ one's have to fit in the first row.

\subsection{Degree of Artin L-functions} It seems that a proportion of $L(\rho,s)$ could be trivial, i.e. $L(\rho,s) = 1$. In order to estimate $N_f$ as accurately as possible, we need to calculate the exact degree of $L(\rho,s)$ for all $\rho$.

Recall that the Artin conductor of $\rho$ is defined as follows: For a Galois extension $L/K$ of local fields and their ring of integers $\mc{O}_K \subseteq \mathcal{O}_L$, the Galois group $G:=$Gal$(L/K)$ has a filtration of higher ramification groups in lower numbering $G_i$. $$G_{-1} = G \supseteq G_0 \supseteq G_1 \supseteq \ldots$$ For the definition of $G_i$ one can see \cite[Ch. VI, pp. 97-106]{serre2013local}. $G_0$ is the inertia group of $L/K$ and $G_1$ is the wild inertia group. Denote $g_i := \lvert G_i\rvert$. For a representation $\rho: G \to GL(V)$, define $$f(\rho, L/K):=\sum_{i \geq 0} \frac{g_i}{g_0}(\dim V - \dim V^{G_i}).$$ $V^{G_i}$ is the invariant subspace of $V$ under $G_i$. Since $G_i = \{1\}$ for $i$ sufficiently large, $f(\rho,L/K)$ is well-defined.

Now for $M/\Ff{q}(t)$ a Galois extension of global function fields (with Galois group $G$ and representation $\rho$), suppose it is ramified at a finite set of primes $S$. Pick any $P \in S$ and fix a prime $\mk{P}$ above $P$, and consider the extension of local fields $M_\mk{P}/\Ff{q}(t)_{P}$, whose Galois group $G_P$ is isomorphic to the decomposition group $D(\mk{P}/P) \subseteq G$. Define $$f(\rho, P) := f(\rho \rvert_{G_P}, M_\mk{P}/\Ff{q}(t)_{P}).$$ Define the Artin conductor of $\rho$
$$f(\rho) :=\sum_{P \in S} f(\rho,P)P.$$ a divisor in $\Ff{q}(t)$.

Let us now return to the original application context. Now $M$ is the splitting field of $f(x)-t$ over $\Ff{q}(t)$, and $x$ is a root of $f(x)-t$. 

The discriminant of field $\Ff{q}(x)/\Ff{q}(t)$ is equal to the discrimiant of $f(x)-t$. The finite ramified primes exactly correspond to the irreducible factors of Disc$(f(x)-t)$, while the prime at infinity $\infty$ is totally ramified in $\Ff{q}(x)$: there is a unique prime ideal $\infty_M$ above the prime at infinity $\infty$ in $\Ff{q}(t)$. For more detailed material discussing the ramification of genus 0 covers of $\mathbb{P}^1$, we refer to \cite{fried1993schur}.


A direct calculation (or an application of Abhyankar's lemma in global function fields) shows: if $P$ is tamely ramified in $\Ff{q}(x)$, then so is it in $M$. When $P$ is tamely ramified in $\Ff{q}(x)$, the expression for local conductor at the finite primes is simplified: $$f(\rho,P) = \dim(V/V^{I(\mk{P}/P)}).$$ for any prime $\mk{P}$ above $P$.

We shall mention here the simplified formula for $\deg L(\rho,s)$ when $\rho$ is tamely ramified everywhere:
$$\deg(L(\rho,s)) = -\dim V + \sum_{\substack{P \text{ finite}\\P \text{ ramified}}} \deg P \cdot \dim(V/V^{I(\mk{P}/P)}).$$

For example, let $V$ be the standard representation of $S_d$, acting on the subspace $y_1+y_2+\ldots+y_d =0$ of complex vector space $\Cc^d$. The inertia subgroup $I(\mk{P}/P)$ is the subgroup of $S_d$ consisting of elements that only interchange the roots $x_1,x_2,\ldots,x_d$ of $f(x)-t$ which map to the same elements in the residue field of $\mk{P}$. Let $W$ be the intersection of family of subspaces $\{H_{i,j}\}$ of $\Cc^d$, where $H_{i,j}$ is defined as such: If $x_i=x_j \mod \mk{P}$, define $H_{i,j} = \{y_i=y_j\}$. Otherwise define $H_{i,j} = \Cc^d$. Then
$$V^{I(\mk{P}/P)} =W \cap V$$ $$\dim(V/V^{I(\mk{P}/P)}) = \sum_{i=1}^g (e_i - 1).$$
$e_i$ is also the ramification index of various primes in $\Ff{q}(x)$ above $P$.
The formula then simplifies to: 
\begin{align*}
    \deg(L(\rho,s)) &= -\dim V + \sum_{\substack{P \text{ finite}\\P \text{ ramified}}} \deg P \cdot \sum_{i}(e_i-1) \\ &=-2\dim V + \sum_{\substack{P \text{ finite}\\P \text{ ramified}}} \deg P \cdot \sum_{i}(e_i-1) + \deg \infty \cdot (d-1) \\ &= -2 \dim V + 2(d-1)=0.
\end{align*}

$e_{\mk{P}}$ is the ramification index of $\mk{P}/P$ and the second to last equality comes from the Riemann-Hurwitz formula (tame ramification case) applied to $\Ff{q}(x)/\Ff{q}(t)$. We have $\deg L(\rho,s) =0$.

We can also study the $L$-function of sign representation $S_n \to \{\pm 1\} \hookrightarrow \Cc^*$. Since $I(\mk{P}/P)$ consists of transpositions of the roots of $f(x)-t$, if it's non-trivial, we have $$\dim(V/V^{I(\mk{P}/P)}) = 1.$$ Hence $$\deg(L(\text{sgn},s)) = -1+\sum_{\substack{P \text{ finite}\\P \text{ ramified}}} \deg P = d-2.$$ Note that the discriminant of $f(x)-t$ is a polynomial of degree $d-1$ in $t$, hence $\sum_{P} \deg P = d-1$.

\begin{remark}
We leave some questions and meaningful remarks.
\begin{itemize}
    \item A natural and important question would be to study the family of functions $L(\rho,s)$ in further detail, e.g. to determine the distribution of $\deg L(\rho,s)$ as $\rho$ ranges over irreducible representations of $S_n$. A nice description of the distribution will allow for more powerful estimates in Theorem \ref{mainthm}.
    \item Let $S$ denote the finite set of places ramified in $\Ff{q}(x)/\Ff{q}(t)$, and let $\mc{F}_\rho$ be the \'etale local system corresponding to $\rho$ on $X_0:=\mathbb{A}^1 \setminus S$. Let $j:X_0 \hookrightarrow \mathbb{P}^1$ be the inclusion morphism. The fact $\deg L(\rho,s) = 0$ is equivalent to $H^i_{et}(\mathbb{P}^1,j_*\mc{F}_\rho) = 0, i =0,1,2$. Notably, this case is extensively studied by Gross \cite{gross2013trivial}.
\end{itemize}
\end{remark}

\section{Application to Quartic Polynomials}

\subsection{Quartic polynomials over finite fields}

Let $f(x)$ be a \textit{generic} polynomial over $\Ff{q}$ of degree 4. To avoid any wild ramification, we require char $\Ff{q} \neq 2,3$. Without loss of generality, suppose $f(0)=0$ and $f$ has been put into the depressed quartic form: $$f(x) = x^4+ax^2+bx.$$ When $b = 0$, then $f$ is essentially a quadratic polynomial in $x^2$, and $f$ cannot be generic. So we need only to focus on $b \neq 0$, in which case $N_f = 5/8q + O(q^{1/2})$.

Let the $S_4$-representations $P_2,P_3$ be as defined in Section 4.1. There are five irreducible representations of $S_4$ in total:

\begin{itemize}
    \item the trivial representation triv, with dimension 1.
    \item the sign representation sgn, with dimension 1.
    \item the standard representation std, with dimension 3.
    \item $R_1:=\text{std}\otimes \text{sgn}$, with dimension 3.
    \item $R_2$ the unique 2-dimensional irreducible representation.
\end{itemize}

We compute the Kostka numbers, and the decomposition of $P_r$ into irreducible representations is as follows:

\begin{itemize}
    \item $P_2$ is of dimension 12, and decomposes into: $$P_2 = \text{triv} \oplus 2 \text{std} \oplus R_1 \oplus  R_2.$$
    \item $P_3 =P_4= \Cc[S_4] = \text{triv} \oplus \text{sgn} \oplus 3 \text{std} \oplus 3 R_1 \oplus 2R_2.$
\end{itemize}

The curves defined in Section 2, $T_2$ and $T_3$ and $T_4$ are, respectively: $$T_2:\frac{x^4-y^4+a(x^2-y^2)}{x-y} = (x^2+y^2+a)(x+y)=-b.$$
an absolutely irreducible plane curve. For general choices of $a$ and $b$ (in fact, when $b^2\neq-8/27\cdot a^3$), it's a smooth plane cubic with geometric genus 1, and its $L$-polynomial has degree $2$. $$T_3:(x^2+y^2+a)(x+y)=-b,x^2+y^2+z^2 +xy+yz+xz=-a.$$ For general choices of $a$ and $b$, it is smooth. Recall that a set-theoretic complete intersection is a variety defined by exactly as many polynomial equations as its codimension. Since smooth set-theoretic complete intersections are also (ideal-theoretic) complete intersections, $T_3$ is an absolutely irreducible smooth curve of geometric genus 4 (complete intersection of multi-degree $(2,3)$), with a degree 8 $L$-polynomial. Hence the degree of Artin $L$-polynomials: $$\sum_{\rho \neq \text{triv}} m_{\rho,2} \deg(L(\rho,s)) = 2.$$ add up to be $2$. $$\sum_{\rho \neq \text{triv}} m_{\rho,3} \deg(L(\rho,s)) = 8.$$
$T_4$ is isomorphic to $T_3$, it is defined in 4-dimensional space: $$T_4: (x^2+y^2+a)(x+y)=-b,x^2+y^2+z^2 +xy+yz+xz=-a,x+y+z+w=0.$$ Note that when $\rho$ is standard representation, the degree of $L(\rho,s)$ is 0. Also, $\deg(L(\text{sgn},s)) = 2$. To ensure the identities above hold, there is only one possibility:

\begin{itemize}
    \item $\deg(L(R_1,s))= 2,\deg(L(R_2,s))= 0, \deg(L(\text{sgn},s))=2$. Theorem \ref{mainthm} shows the implied constant is: \begin{align*}
      \lvert m_{R_1,2}/2! -m_{R_1,3}/3!+m_{R_1,4}/4! \rvert \cdot \deg L(R_1,s) + & \lvert m_{\text{sgn},2}/2! -m_{\text{sgn},3}/3!+m_{\text{sgn},4}/4! \rvert \cdot \deg L(\text{sgn},s) \\ &= (1/2! - 3/3!+3/4!) \cdot 2 + (1/3!-1/4!) \cdot 2 \\ &= 1/2.
    \end{align*} So \begin{align*}
        \lvert N_f - \frac{5}{8}q \vert & \leq \frac{1}{2}\sqrt{q} + O_d(1).
    \end{align*} The implied constant in $O(\sqrt{q})$ matches our experimental evidence very well. We need to compute $O_d(1)$ precisely. Recall the formula from Section 2: \begin{align*}
        N_f = q + \sum_{r=2}^d (-1)^{r+1} \frac{\#C_r(\Ff{q})}{r!} + \sum_{r=2}^d (-1)^{r+1} \frac{n_r'-\#C_r(\Ff{q})}{r!}.
    \end{align*} And the error term of Theorem \ref{mainthm} is precisely $$\sum_{r=2}^4 (-1)^{r+1} \frac{n_r'-\#C_r(\Ff{q})}{r!}+(\mu_4-1) = \sum_{r=2}^4 (-1)^{r+1} \frac{(n_r'- n_r) + (n_r -\#C_r(\Ff{q}))}{r!}-\frac{3}{8}.$$ Since $T_2$, $T_3$ should have no singularities, the term $n_r - \#C_r(\Ff{q})$ simply comes from the points at infinity of $T_2$, $T_3$ and $T_4$. In our case, $T_2$ has at most 3 points at infinity, and $T_3$ and $T_4$ has at most 6 points at infinity (Note that they have the same number of points at infinity). So $$\lvert\sum_{r=2}^4 (-1)^{r+1} \frac{n_r - \#C_r(\Ff{q})}{r!} \rvert= 3/2! - a/3!+a/4! \leq 3/2.$$ where $a$ is the number of infinity points of $T_3$. the term $n'_r-n_r$ on the other hand, comes from the intersection of hyperplanes $x=y$, $y=z$, etc. with the curves $T_2$, $T_3$ and $T_4$. For $T_2$, we need to rule out the points of $T_r$ on $x=y$ ($\leq 3$ points). Similarly, for $T_3$ and $T_4$ we need to rule out points on 3 hyperplanes and 6 hyperplanes, respectively (on each hyperplane, we rule out $\leq 6$ points). To sum up: $$\lvert \sum_{r=2}^4 (-1)^{r+1} \frac{n_r'-n_r}{r!} \rvert = \lvert a_1/2 - a_2/3! + a_3/4! \rvert \leq 3.$$ where $0 \leq a_1 \leq 3, 0 \leq a_2 \leq 3 \times 6= 18, 0 \leq a_3 \leq 36$.
\end{itemize}

This concludes our proof of Corollary \ref{quar}.

\begin{remark}
    For quartic polynomials, our results are rather impressive. However, when we attempt to apply the same method to quintic polynomials, the situation becomes more difficult: the irreducible representations of $S_5$ are more complicated and the decomposition of $P_r$ is much less manageable, making it hard to obtain further improvements. We leave the case of quintic polynomial to the readers.
\end{remark}

\bibliography{biblio}

\end{document}